\documentclass[11pt, leqno]{article}
\usepackage{euscript,amstext,amsmath,amsthm,a4,amssymb}

\newcommand{\comment}[1]{}

\newtheorem{thm}{Theorem}[section]
\newtheorem{prop}[thm]{Proposition}
 \newtheorem{cor}[thm]{Corollary}
\newtheorem{lemma}[thm]{Lemma}
\theoremstyle{remark}
\newtheorem{Rem}[thm]{Remark}
\theoremstyle{definition}

\title{Real embeddings, $\eta$-invariant\\  and  Chern-Simons current\thanks{Partially
supported by the Ministry  of Education and the National Natural
Science Foundation of China.}}
\author{Huitao Feng\footnote{Chern
Institute of Mathematics \& LPMC, Nankai University, Tianjin
300071, P.R. China. (fht@nankai.edu.cn)},\ \ \ Guangbo
Xu\footnote{Chern Institute of Mathematics \& LPMC, Nankai
University, Tianjin 300071, P.R. China. (0210042@nankai.edu.cn)}\
\ \ and\ \ Weiping Zhang\footnote{Chern Institute of Mathematics
\& LPMC, Nankai University, Tianjin 300071, P.R. China.
(weiping@nankai.edu.cn)}}
\date{}

\begin{document}

\maketitle
\begin{abstract} We present an alternate proof of the
Bismut-Zhang localization  formula for $\eta$-invariants without
using  the  analytic   techniques developed by Bismut-Lebeau. A
Riemann-Roch property for Chern-Simons currents, which is of
independent interest,  is established in due course.

\comment{$\ $


\noindent {\bf 2000 MR Subject Classification} 58J}
\end{abstract}

\renewcommand{\theequation}{\thesection.\arabic{equation}}
\setcounter{equation}{0}

\section{Introduction} \label{s1}

The $\eta$ invariant of Atiyah-Patodi-Singer was introduced in
\cite{APS} as the correction term on the boundary of the index
theorem for Dirac operators on manifolds with boundary. Since then
it has appeared in many parts of geometry, topology as well as
physics. We first recall the definition of this important
invariant.

Let $M$ be an odd dimensional oriented closed  spin manifold
carrying a Riemannian metric $g^{TM}$. Let $S(TM)$ be the
associated Hermitian bundle of spinors. Let $\mu$ be a Hermitian
vector bundle over $M$ carrying a unitary connection.
 Let
\begin{align}\label{1.1} D^{\mu }:\Gamma(S(TM)\otimes \mu )\longrightarrow
\Gamma(S(TM)\otimes \mu )
\end{align}
denote the corresponding (twisted) Dirac operator, which is
formally self-adjoint (cf. \cite{BGV}).

For any $s\in {\bf C}$ with ${\rm Re}(s)>>0$, following
\cite{APS}, set
\begin{align}\label{1.2}
\eta\left(D^{\mu  },s\right)=\sum_{\lambda\in{\rm Spec}(D^{\mu
})\setminus \{ 0\}}{{\rm Sgn}(\lambda)\over |\lambda|^s}.
\end{align} By \cite{APS}, one knows that $\eta (D^{\mu  },s
)$ is a holomorphic function in $s$ when ${\rm Re}(s)>{\dim M\over
2}$. Moreover, it extends to  a meromorphic function over ${\bf
C}$, which is holomorphic at $s=0$. The $\eta$ invariant of
$D^{\mu }$, in the sense of Atiyah-Patodi-Singer \cite{APS}, is
defined by
\begin{align}\label{1.3}\eta\left(D^{\mu } \right)=
\eta\left(D^{\mu  },0\right) ,
\end{align}
while the corresponding {\it reduced} $\eta$ invariant is defined
and denoted  by
\begin{align}\label{1.4}\overline{\eta}\left(D^{\mu  } \right)=
{\dim \left(\ker D^{\mu  }\right)+\eta\left(D^{\mu  }\right)\over
2} .
\end{align}

Let $i:Y\hookrightarrow X$ be an embedding between two odd
dimensional compact oriented spin Riemannian manifolds. For any
Hermitian vector bundle $\mu$ over $Y$ carrying a Hermitian
connection, Bismut and Zhang \cite[Theorem 2.2]{BZ} established a
mod {\bf Z} formula, expressing $\overline{\eta}(D^\mu)$ through
the $\eta$-invariants associated to certain   direct image
$i_!\mu$ in the sense of Atiyah-Hirzebruch \cite{AH}, up to some
geometric Chern-Simons   current. This formula, in some sense,
might be thought of as a Riemann-Roch type formula for
$\eta$-invariants under embeddings.

The proof in \cite{BZ} relies  heavily on the analytic techniques
developed in the difficult paper of Bismut and Lebeau \cite{BL}.
On the other hand, in a special case where $X$ is certain higher
dimensional sphere, a more geometric proof of the above
Bismut-Zhang formula was given in \cite{Z2} by making use of the
mod $k$ index theorem of Freed-Melrose \cite{FM}. As a
consequence, one gets a purely geometric formula for the mod {\bf
Z} part of $\overline{\eta}(D^\mu)$ (cf. \cite[Theorem 2.2]{Z2}).

It is natural to ask   whether the original Bismut-Zhang
localization  formula for arbitrary $X$  can also be proved
without using the techniques developed in \cite{BL}.  The purpose
of this paper is to show that this is indeed the case. More
precisely, as indicated in \cite[Remark 3.2]{Z2}, we will embed
$X$ into a sufficiently high dimensional odd sphere $S^{2N-1}$ and
apply the proved case to $Y\hookrightarrow S^{2N-1}$ and
$X\hookrightarrow S^{2N-1}$ respectively, to get the final
formula. Meanwhile, we also establish  a Riemann-Roch type formula
for the involved Chern-Simons currents (cf. Theorem \ref{t3.4}),
which has  its own interest.

The rest of this paper is organized as follows. In Section 2, we
recall the geometric construction  of the direct image $i_!\mu$
and the Bismut-Zhang localization formula. In Section 3, we
present our alternate proof of the Bismut-Zhang  formula.

\section{Direct image and the Bismut-Zhang localization formula for $\eta$-invariants}\label{s2}

This section is organized as follows. In Section \ref{2a}, we
recall some basic notions of super vector bundles. In Section
\ref{2b}, we recall the geometric construction of the direct image
of a vector bundle under embeddings, as well as the associated
Chern-Simons current.  In Section \ref{2c}, we recall the statement of the Bismut-Zhang localization formula for
$\eta$-invariants.

\subsection{Basic notions of super vector bundles}\label{2a} \setcounter{equation}{0}

Let $\xi=\xi_+\oplus \xi_-$ be a ${\bf Z}_2$-graded Hermitian
vector bundle in the sense of Quillen \cite{Q} over a manifold.

Let $v:\xi_+\rightarrow \xi_-$ be an endomorphism between the
vector bundles $\xi_+$ and $\xi_-$. Let $v^*:\xi_-\rightarrow
\xi_+$ be the adjoint of $v$ with respect to the Hermitian metrics
on $\xi_\pm$. Then $V=v+v^*:\xi\rightarrow \xi$ is an odd self-adjoint
endomorphism of the Hermitian super vector bundle $\xi$. We will
denote this set of data by $(\xi_+,\xi_-,V)$.

Let  $ {\rm Supp}(V)$ denote the subset where $V$ is not
invertible.

For two super vector bundles with  odd endomorphism
$((\xi_1)_+,(\xi_1)_-,V_1)$ and
$(\left(\xi_2\right)_+,(\xi_2)_-,V_2)$, we can form their direct
sum
\begin{align}\label{2.1}\left((\xi_1)_+,(\xi_1)_-,V_1\right)\oplus
\left((\xi_2)_+,(\xi_2)_-,V_2\right)=\left((\xi_1)_+\oplus(\xi_2)_+,(\xi_1)_-\oplus(\xi_2)_-,
V_1\oplus V_2\right)\end{align}
 and also the \emph{super} tensor
product
\begin{multline}\label{2.2}
\left(\left(\xi_1\right)_+,\left(\xi_1\right)_-,V_1\right)
\widehat{\otimes}\left(\left(\xi_2\right)_+,\left(\xi_2\right)_-,V_2\right)\\
=\left(\left(\xi_1\right)_+\otimes\left(\xi_2\right)_+\oplus\left(\xi_1\right)_-
\otimes\left(\xi_2\right)_-,\left(\xi_1\right)_+\otimes\left(\xi_2\right)_-
\oplus\left(\xi_1\right)_-\otimes\left(\xi_2\right)_+, V_1\otimes
{\rm Id}_{\xi_2}+{\rm Id}_{\xi_1}\otimes V_2\right),
\end{multline} with obvious induced Hermitian metrics.

The following formulas are clear from the definitions,
\begin{align}\label{2.3}\mathrm{Supp}\left(V_1\oplus V_2\right)=\mathrm{Supp}\left(V_1\right)\cup
\mathrm{Supp}\left(V_2\right),\end{align} while
\begin{align}\label{2.4}\mathrm{Supp}\left(V_1\otimes
{\rm Id}_{\xi_2}+{\rm Id}_{\xi_1}\otimes
V_2\right)=\mathrm{Supp}\left(V_1\right)\cap\mathrm{Supp}
\left(V_2\right).\end{align}

\subsection{Geometric construction of direct images and the
associated Chern-Simons current}\label{2b}

For completeness of this paper, we recall the geometric
constructions of direct images and the associated Chern-Simons
currents from \cite{BZ} and \cite{Z2}.

Let $i:Y\hookrightarrow X$ be an embedding between two closed
oriented  spin manifolds. We make the assumption that $\dim X-\dim
Y$ is even and that if $N$ denotes the normal bundle to $Y$ in
$X$, then $N$  is orientable,  spin and carries an induced
orientation as well as a (fixed) spin structure.

Let $g^N$ be a Euclidean metric on $N$ and $\nabla^N$  a Euclidean
connection on $N$ preserving $g^N$. Let $S(N)$ be the vector
bundle of spinors associated to $(N,g^N)$. Then $S(N)=S_+(N)\oplus
S_-(N)$ (resp. its dual $S^*(N)=S^*_+(N)\oplus S^*_-(N)$) is a
${\bf Z}_2$-graded complex vector bundle over $Y$ carrying   an
induced Hermitian metric $g^{S(N)}=g^{S_+(N)}\oplus g^{S_-(N)}$
(resp. $g^{S^*(N)}=g^{S^*_+(N)}\oplus g^{S^*_-(N)}$) from $g^N$,
as well as a Hermitian connection
$\nabla^{S(N)}=\nabla^{S_+(N)}\oplus\nabla^{S_-(N)}$ (resp.
$\nabla^{S^*(N)}=\nabla^{S^*_+(N)}\oplus\nabla^{S^*_-(N)}$)
induced from $\nabla^N$.

For any $r>0$, set $N_r=\{Z\in N: |Z|< r\}.$ Then there is
$\varepsilon_0>0$ such that $N_{2\varepsilon_0}$ is diffeomorphic
to an open neighborhood of $Y$ in $X$. Without confusion we now
view directly $N_{2\varepsilon_0}$ as an open neighborhood of $Y$
in $X$.

Let $\pi:N\rightarrow Y$ denote the projection of the normal
bundle $N$ over $Y$.

If $Z\in  {N}$, let $\tilde{c}(Z)\in {\rm End}(S^*({N})) $ be the
transpose of $c(Z)$ acting on $S( {N})$.

Let $\tau^{ {N}*}\in {\rm End}(S^*( {N})) $ be the transpose of
$\tau^{ {N}}$ defining the ${\bf Z}_2$-grading of $S( {N})=S_+(
{N})\oplus S_-( {N})$.

Let $\pi^*(S^*(N))$ be the pull back bundle of $S^*(N)$ over $N$.

For any $Z\in N$ with $Z\neq 0$, let $\tau^{ {N}*}\tilde{c}(Z):
\pi^*(S^*_\pm(N))|_Z\rightarrow \pi^*(S^*_\mp(N))|_Z$ denote the
corresponding pull back isomorphisms at $Z$.

Let  $(\mu,g^\mu)$ be a Hermitian vector bundle over $Y$ carrying
a Hermitian connection $\nabla^\mu$.

In this paper, by a direct image of $\mu$ under the embedding
$i:Y\to X$, we always mean a triple $(\xi_+,\xi_-,V)$ described in
Section 2.1 with a Hermitian connection
$\nabla^\xi=\nabla^{\xi_+}\oplus\nabla^{\xi_-}$ verifying the
following fundamental assumptions (cf. \cite[(1.10)-(1.12)]{BZ}):

1) $V$ is invertible on $X\setminus Y$ and $(\ker V)|_Y$ has a
constant dimension;

2) the following identification
\begin{align}\label{asp1}
\left(\pi^*(\ker V)|_Y,\pi^*g^{(\ker
V)|_Y},\dot{\partial}_ZV\right)\simeq \left(\pi^*\left(\mu\otimes
S^*(N)\right),\pi^*g^{\mu\otimes
S^*(N)},\tau^{N*}\tilde{c}(Z)\right)
\end{align}
holds over  $ N $, where the map $\dot{\partial}_ZV$ is defined by
\begin{align}\label{asp2}
\dot{\partial}_ZV=P^{\ker V}\left({\partial}_ZV\right)P^{\ker V}
\end{align}with
respect to any smooth trivialization of $\xi$ near $\pi(Z)$ and
$P^{\rm ker V}$ denotes the orthogonal projection from $\xi$ onto
${\ker V}$;

3) under the identification (2.5) the following connections
identification holds,
\begin{align}\label{asp3}
\nabla^{(\ker V)|_Y}=\nabla^{\mu\otimes S^*(N)},
\end{align}
where $\nabla^{(\ker V)|_Y}$ is defined by
\begin{align}\label{asp3}
\nabla^{(\ker V)|_Y}=P^{\ker V}\nabla^{\xi|_Y}P^{\ker V}.
\end{align}

Clearly, $\xi_+-\xi_-\in \widetilde{K}(X)$ is exactly the
Atiyah-Hirzebruch direct image $i_!\mu$ of $\mu$ constructed in [2].

\bigskip
We now describe a concrete  realization of the direct image  of
$\mu$ for an embedding $i:Y\to X$ which verifies the assumption
1)--3) (see also   \cite{Z2}).

Let $(F,g^F)$ be a Hermitian vector bundle over $Y$ carrying a
Hermitian connection $\nabla^F$ such that $S_-(N)\otimes\mu\oplus F$
is a trivial complex vector bundle over $Y$ (cf. \cite{A}). Then
\begin{align}\label{2.9}\tau^{ {N}*}\tilde{c}(Z)\oplus\pi^*{\rm Id}_F
:\pi^*\left(S_+^*(N)\otimes\mu\oplus F\right)\rightarrow
\pi^*\left(S_-^*(N)\otimes\mu\oplus F\right)\end{align} induces an
isomorphism between two trivial  vector bundles over $
N_{2\varepsilon_0}\setminus Y$.

Clearly, $\pi^*(S_\pm^* (N)\otimes\mu\oplus F)|_{\partial
N_{2\varepsilon_0}}$ extend smoothly to two trivial complex vector
bundles over $X\setminus N_{2\varepsilon_0}$. Moreover, the
isomorphism $\tau^{ {N}*}\tilde{c}(Z)\oplus\pi^*{\rm Id}_F$ over
$\partial N_{2\varepsilon_0}$ extends smoothly to an isomorphism
between these two trivial vector bundles over $X\setminus
N_{2\varepsilon_0}$.

In summary, what we get is a ${\bf Z}_2$-graded Hermitian vector
bundle $(\xi=\xi_+\oplus \xi_-, g^\xi=g^{\xi_+}\oplus g^{\xi_-})$
over $X$ such that
\begin{align}\label{2.10}\xi_\pm|_{N_{\varepsilon_0}}
=\pi^*\left.\left(S_\pm^*(N)\otimes\mu\oplus
F\right)\right|_{N_{\varepsilon_0}},\ \ \ \
g^{\xi_\pm|_{N_{\varepsilon_0}}}
=\pi^*\left.\left(g^{S_\pm^*(N)\otimes\mu}\oplus
g^F\right)\right|_{N_{\varepsilon_0}},\end{align} where
$g^{S_\pm^*(N)\otimes\mu}$ is the tensor product Hermitian metric
on $S_\pm^*(N)\otimes\mu$ induced from $g^{S_{\pm}^*(N)}$ and
$g^\mu$.

It is easy to see that there exists an odd self-adjoint
automorphism $V$ of $\xi$ such that
\begin{align}\label{2.11}V|_{N_{\varepsilon_0}} =\tau^{ {N}*}\widetilde{c}(Z)\oplus\pi^*{\rm
Id}_F.\end{align}
Moreover, there is a ${\bf Z}_2$-graded
Hermitian connection $\nabla^\xi=\nabla^{\xi_+}\oplus
\nabla^{\xi_-}$ on $\xi=\xi_+\oplus\xi_-$ over $X$ such that
\begin{align}\label{2.12}
\nabla^{\xi_\pm}|_{N_{\varepsilon_0}}=
\pi^*\left(\nabla^{S^*_\pm(N)\otimes\mu}\oplus\nabla^F\right),\end{align}
where $\nabla^{S^*_\pm(N)\otimes\mu}$ is the Hermitian connection
on $\nabla^{S^*_\pm(N)\otimes\mu}$ defined by
$\nabla^{S^*_\pm(N)\otimes\mu}=\nabla^{S^*_\pm(N)}\otimes {\rm
Id}_\mu +{\rm Id}_{S^*_\pm (N)}\otimes\nabla^\mu$.

Clearly, the fundamental assumptions (\ref{asp1}) and (\ref{asp2})
hold for this geometric construction. We will call $(\xi_+,\xi_-,
V)$ constructed as such a {\it geometric} direct image of $\mu$.

\bigskip
Before recalling the construction of  the associated Chern-Simons
current, we give some notation first.

Let $i^{1/2}$ be a fixed square root of $i=\sqrt{-1}$. The objects
which will be considered in the sequel do not depend on this
square root. Let $\varphi$ be the map $\alpha\in
\Lambda^*(T^*X)\rightarrow (2\pi i)^{-{\deg \alpha\over 2}}\alpha
\in\Lambda^*(T^*X)$.

If $E$ is a real vector bundle over $X$ carrying with a connection
$\nabla^E$, we denote by $\widehat{A}(E,\nabla^E)$ the Hirzebruch
characteristic form defined by
\begin{align}\label{2.13}\widehat{A}(E, \nabla^{E})
={\det}^{1/2}\left({{\sqrt{-1}\over  4\pi}R^{E} \over \sinh\left({
\sqrt{-1}\over 4\pi}R^{E}\right)}\right)=\varphi\,
{\det}^{1/2}\left({{ R^{E}\over  2 }  \over \sinh\left({
 {R^{E}}\over 2 } \right)}\right) ,\end{align} where
$R^E=\nabla^{E,2}$ is the curvature of $\nabla^E$. While if $E'$
is a complex vector bundle carrying with a connection
$\nabla^{E'}$, we denote by ${\rm ch}(E',\nabla^{E'})$ the Chern
character form associated to $(E',\nabla^{E'})$ (cf. \cite[Section
1]{Z1}).

For $T\geq 0$, let $C_T$ be the superconnection on the super
vector bundle $\xi$ defined by
\begin{align}\label{2.14}C_T=\nabla^\xi+\sqrt{T}V.\end{align}
The curvature $C_T^2$ of $C_T$ is a smooth section of
$(\Lambda^*(T^*X)\widehat{\otimes}{\rm End}(\xi))^{\rm even}$.

By \cite{Q}, we know that for any $T>0$,
\begin{align}\label{2.15}{\partial\over\partial T}{\rm
Tr}_s\left[\exp\left(-C^2_T\right)\right]=-{d\over 2\sqrt{T}}{\rm
Tr}_s\left[V\exp\left(-C^2_T\right)\right].\end{align}

 By proceeding  as
in \cite{B1}, \cite{B2} and \cite[Definition 1.3]{BZ}, one can
construct the Chern-Simons current $\gamma^{\xi, V}$ as
\begin{align}\label{2.16}\gamma^{\xi,
V}={1\over \sqrt{2\pi i}}\int_0^{+\infty}\varphi{\rm
Tr}_s\left[V\exp\left(-C^2_T\right)\right]{dT\over
2\sqrt{T}}.\end{align}

Let $\delta_Y$ denote the current of integration over the oriented
submanifold $Y$ of $X$.

By \cite[Theorem 1.4]{BZ}, we have that
\begin{align}\label{2.17}d\gamma^{\xi,
V}={\rm ch}\left(\xi_+,\nabla^{\xi_+}\right)-{\rm
ch}\left(\xi_-,\nabla^{\xi_-}\right)-{\widehat{A}}^{-1}\left(N,\nabla^N\right){\rm
ch }\left(\mu,\nabla^\mu\right)\delta_Y.\end{align} Moreover, as
indicated in \cite[Remark 1.5]{BZ}, by proceeding as in
\cite[Theorem 3.3]{BGS}, one can prove that $\gamma^{\xi, V}$ is a
locally integrable current.

\subsection{The Bismut-Zhang localization formula for $\eta$
invariants}\label{2c}

We assume in this subsection that $i:Y\hookrightarrow X$ is an
embedding between two odd dimensional closed oriented spin
manifolds. Then the normal bundle $N$ to $Y$ in $X$ is even
dimensional and carries a canonically induced orientation and spin
structure.

Let $g^{TX}$  be a Riemannian metric on $TX$. Let $g^{TY}$ be the
restricted Riemannian metric on $TY$. Let $\nabla^{TX}$ (resp.
$\nabla^{TY}$) denote the Levi-Civita connection associated to
$g^{TX}$ (resp. ${g^{TY}}$). Without loss of generality we may and
we will make the assumption that the embedding
$(Y,g^{TY})\hookrightarrow (X,g^{TX})$ is totally geodesic.

Let $N$ carry the canonically induced Euclidean metric as well as
the Euclidean connection.

The definition of the reduced $\eta$ invariant for a
(twisted) Dirac operator on an odd dimensional spin Riemannian
manifold has been recalled in Section 1.

Under our assumptions, we can state the Bismut-Zhang localization
formula for $\eta$-invariants \cite{BZ}   as follows, of which a
special case was proved in \cite{B2}.

\begin{thm}\label{t2.1}{\rm (Bismut-Zhang \cite[Theorem 2.2]{BZ})}
If $(\xi_+,\xi_-,V)$ is a direct
image of $\mu$ for a totally geodesic embedding $i:Y\hookrightarrow X$, then the following identity holds,
\begin{align}\label{2.18}\overline{\eta}\left(D^{\xi_+}\right) -
\overline{\eta}\left(D^{\xi_-}\right)
\equiv\overline{\eta}\left(D^\mu\right)+\int_X\widehat{A}\left(TX,\nabla^{TX}\right)\gamma^{\xi,V}\
\ \ {\rm mod}\ {\bf Z}.\end{align}\end{thm}

\begin{Rem}\label{t2.2} The extra Chern-Simons form in \cite[Theorem 2.2]{BZ} disappears
here simply because we have made the simplifying assumption that
the isometric embedding $(Y,g^{TY})\hookrightarrow (X,g^{TX})$ is
totally geodesic.\end{Rem}

\begin{Rem}\label{t2.3} The proof of (\ref{2.18}) in \cite{BZ} relies heavily on the analytic techniques developed
in a difficult paper of Bismut-Lebeau \cite{BL}.  By using the mod
$k$ index theorem of Freed-Melrose \cite{FM}, Zhang  showed in
\cite{Z2}
 that there exists an embedding $i:Y\hookrightarrow S^{2m-1}$
of $Y$ to some higher dimensional sphere and a geometric direct
image $(\xi^\prime_+,\xi^\prime_-,V^\prime)$ of $\mu$ on
$S^{2m-1}$ such that the    following special case of (\ref{2.18})
holds,
\begin{align}\label{2.19}\overline{\eta}\left(D^{\xi^\prime_+}\right) -
\overline{\eta}\left(D^{\xi^\prime_-}\right)
\equiv\overline{\eta}\left(D^\mu\right)+\int_{S^{2m-1}}
\widehat{A}\left(TS^{2m-1},\nabla^{TS^{2m-1}}\right)\gamma^{\xi^\prime,V^\prime}\
\ \ {\rm mod}\ {\bf Z},\end{align} without using the techniques of
Bismut-Lebeau \cite{BL}.

In the rest of this paper, we will present a proof of (\ref{2.18})
by using (\ref{2.19}), which still avoids the use of the
techniques in \cite{BL}.
 \end{Rem}

\section{A proof of Theorem \ref{t2.1}}\label{s3}\setcounter{equation}{0}

In this section we will present an alternate proof of Theorem
\ref{t2.1} by embedding $X$ to a higher dimensional sphere.

\bigskip
Let $(\xi_+,\xi_-,V)$ be a direct image of $\mu$ for an embedding
$i:Y\hookrightarrow X$ in the sense of Section \ref{2b}.

For the sake of convenience, we denote by $H_X(\xi_+,\xi_-,V)\in
{\bf R}/{\bf Z}$ the quantity defined by
\begin{align}\label{3.1} H_X\left(\xi_+,\xi_-,V\right)\equiv
\overline{\eta}\left(D^{\xi_+}\right)-
\overline{\eta}\left(D^{\xi_-}\right)-\int_X\widehat{A}\left(TX,\nabla^{TX}\right)\gamma^{\xi,V}-
\overline{\eta}\left(D^\mu\right)\ \ {\rm mod}\ {\bf Z}.
\end{align}
Then (\ref{2.18}) and (\ref{2.19}) can be rewritten as
\begin{align}\label{3.2}
  H_X\left(\xi_+,\xi_-,V\right)=0,
\end{align}
and
\begin{align}\label{3.3}
 H_{S^{2m-1}}\left(\xi^\prime_+,\xi^\prime_-,V^\prime\right)=0
\end{align}
respectively.

The rest of this section is organized as follows. In Section
\ref{3.h}, we   prove two basic properties of
$H_X(\xi_+,\xi_-,V)$. In Section \ref{3a}, we describe a
constructions of direct images under successive embeddings and the
associated Chern-Simons current. In Section \ref{3b}, we study the
relations between the Chern-Simons currents constructed in Section
\ref{3a} and establish a Riemann-Roch type formula for them. In
Section \ref{3c}, we use (\ref{3.3}) and the Riemann-Roch type
formula established in Section \ref{3b}  to give an alternative
proof of the Bismut-Zhang localization formula.

\subsection{Basic properties of the $H$-quantity}\label{3.h}
In this subsection, we will prove two properties of the
$H$-quantities defined above, from which one can deduce that the
$H$-quantity  depends only on the isotropy class of the
embeddings.

\begin{lemma}\label{l3.2} If  $i_s:Y_s\hookrightarrow X_s$, $0\leq
s\leq 1$, is a smooth family of embeddings between odd dimensional
compact oriented spin Riemannian manifolds such that
$i_0:Y_0\hookrightarrow X_1$, $i_1:Y_1\hookrightarrow X_1$ are
totally geodesic embeddings,  and $(\xi_{s,+},\xi_{s,-},V_s)$ is a
smooth family of   direct images of the complex vector bundle
$\mu_s$ over $Y_s$, then the following identity in ${\bf R}/{\bf
Z}$ holds,
\begin{align}\label{3.19}
H_X(\xi_{0, +},\xi_{0, -},V_0)= H_X(\xi_{1, +},\xi_{1, -},V_1) .
\end{align}
\end{lemma}
{\it Proof.} Without loss of generality, we may and we will assume
that the above smooth family is a locally  constant family for
$s\in [0,{1\over 8}]\cup [{7\over 8},1]$.

Set $Y=Y_0$ and ${\widehat Y}=I\times Y$.

We equip $\widehat Y$ with the metric $ds^2\oplus g^{TY_s}$,
$0\leq s\leq 1$. Let $\nabla^{\widehat Y}$ denote the associated
Levi-Civita connection.

Clearly, $\widehat Y$ is an oriented spin even-dimensional
Riemannian manifold with boundary $\partial\widehat Y={\overline
Y}_1\bigcup Y_0$, where ${\overline Y}_1$ is a copy of $Y_1$ but
with the reversed orientation.

Let ${\widehat\mu}$ be the canonical complex vector bundle over
${\widehat Y}$ such that ${\widehat\mu}|_{\{s\}\times Y}=\mu_s$,
let $g^{\widehat\mu}$ (resp. $\nabla^{\widehat\mu}$) be the
Hermitian metric (resp. connection) on ${\widehat\mu}$ such that
$g^{\widehat\mu}|_{\mu_s}=g^{\mu_s}$ (resp.
$\nabla^{\widehat\mu}|_{\mu_s}=\nabla^{\mu_s}$).

By the Atiyah-Patodi-Singer index  theorem \cite{APS}, one has
\begin{align}\label{3.20}
\int_{\widehat Y}{\widehat A}\left(T{\widehat Y},\nabla^{\widehat
Y}\right){\rm ch}\left({\widehat\mu},\nabla^{\widehat\mu}\right)-
\overline\eta\left(D^{\widehat\mu}_{\partial
\widehat{Y}}\right)\equiv 0, \ \mbox{\rm mod $\bf Z$}.\end{align}

Let $X=X_0$ and ${\widehat X}=I\times X$.

It is easy to see that  there exists a metric $g^{T{\widehat X}}$
on $T\widehat{X}$ such that it equals to $ds^2\oplus g^{TX_s}$ on
$[0,{1\over 8}]\times X\cup [{7\over 8},1]\times X$, and that the
canonical embedding $i_{\widehat Y}:{\widehat Y}\hookrightarrow
{\widehat X}$ is totally geodesic. Let $\nabla^{\widehat X}$
denote the associated Levi-Civita connection.

Clearly, $\widehat X$ is an oriented spin even-dimensional
Riemannian manifold with boundary $\partial\widehat X={\overline
X}_1\bigcup X_0$, where ${\overline X}_1$ is a copy of $X_1$ but
with the reversed orientation.

The smooth family of direct images $(\xi_s,V_s)$ also lifts to
$\widehat{X}$ canonically and form a direct image ($\widehat{\xi},
\widehat{V}$) over $\widehat{X}$ of $\widehat{\mu}$ over
${\widehat Y}$, which is of product structure near the boundary.
In particular, when restricted to $s\in[0,{1\over 8}] \cup
[{7\over 8},1]$, $({\widehat\xi}_s,V^{\widehat\xi}_s) $ is the
(locally constant) geometric direct image of $\mu_s\to Y_s$ for
the totally geodesic embedding $i_s:Y_s\to X_s$.

Denote by ${\widehat N}$   the normal bundle to ${\widehat Y}$ in
${\widehat X}$.

By (\ref{2.17}), one has
\begin{align}\label{3.22}
d\gamma^{{\widehat\xi},{\widehat V}}={\rm
ch}\left({\widehat\xi}_+,\nabla^{{\widehat\xi}_+}\right) -{\rm
ch}\left({\widehat\xi}_-,\nabla^{{\widehat\xi}_-}\right)-
{\widehat A}^{-1}\left({\widehat N},\nabla^{\widehat N}\right)
{\rm
ch}\left(\widehat\mu,\nabla^{\widehat\mu}\right)\delta_{\widehat
Y}.\end{align}

On the other hand, by the Atiyah-Patodi-Singer index theorem
\cite{APS}, one has
\begin{align}\label{3.23}
\overline\eta\left(D^{\widehat{\xi}}_{\partial{\widehat X}}\right)
\equiv\int_{\widehat X}{\widehat A}\left(T{\widehat
X},\nabla^{\widehat X}\right) {\rm
ch}\left({\widehat\xi},\nabla^{\widehat\xi}\right)\ \ {\rm mod}\
{\bf Z}.\end{align}

From (\ref{3.20}), (\ref{3.22}) and (\ref{3.23}), one gets
\begin{multline}\label{4.0}
 \overline\eta\left(D^{\xi_{0,+}}\right)-\overline{\eta}\left(D^{\xi_{0,-}}\right)
-\left(\overline\eta\left(D^{\xi_{1,+}}\right)-\overline{\eta}\left(D^{\xi_{1,-}}\right)\right)\\
\equiv\int_{\widehat X}{\widehat A}\left(T{\widehat
X},\nabla^{\widehat X}\right) \left({\rm
ch}\left({\widehat\xi}_+,\nabla^{{\widehat\xi}_+}\right)-
{\rm ch}\left({\widehat\xi}_-,\nabla^{{\widehat\xi}_-}\right)\right)\\
\equiv\int_{\widehat X}{\widehat A}\left(T{\widehat
X},\nabla^{T\widehat X}\right) d\gamma^{{\widehat\xi},{\widehat
V}}+ \int_{\widehat X}{\widehat A}\left(T{\widehat
X},\nabla^{T\widehat X}\right)
{\widehat A}^{-1}\left({\widehat N},\nabla^{\widehat N}\right)
{\rm ch}\left(\widehat\mu,\nabla^{\widehat\mu}\right)\delta_{\widehat Y}\\
\equiv\int_X{\widehat
A}\left(TX,\nabla^{TX}\right)\gamma^{\xi_0,V_0}-\int_X{\widehat
A}\left(TX,\nabla^{TX}\right)\gamma^{\xi_1,V_1}+ \int_{\widehat
Y}{\widehat A}\left(T{\widehat Y},\nabla^{T\widehat Y}\right)
 {\rm ch}\left(\widehat\mu,\nabla^{\widehat\mu}\right) \\
\equiv\int_X{\widehat
A}\left(TX,\nabla^{TX}\right)\gamma^{\xi_0,V_0}-\int_X{\widehat
A}\left(TX,\nabla^{TX}\right)\gamma^{\xi_1,V_1}
+\overline\eta\left(D^{\widehat\mu}_{\partial{\widehat Y}}\right)\\
\equiv\int_X{\widehat
A}\left(TX,\nabla^{TX}\right)\gamma^{\xi_0,V_0}-\int_X{\widehat
A}\left(TX,\nabla^{TX}\right)\gamma^{\xi_1,V_1}+\overline{\eta}
\left(D^{\mu_0}\right)-\overline{\eta}\left(D^{\mu_1}\right)\ \
\mbox{\rm mod}\ {\bf Z}.
\end{multline}
From (\ref{3.1}) and (\ref{4.0}), we get (\ref{3.19}).\ \ Q.E.D.

$\ $

\begin{lemma}\label{l3.1}
Let $i:Y\hookrightarrow X$ be a totally geodesic embedding between
the odd dimensional compact oriented spin Riemannian manifolds and
$ \mu$ a Hermitian vector bundle over $Y$ carrying a Hermitian
connection. Then for any two direct images
$(\xi_{k,+},\xi_{k,-},V_k)$, $k=1,2$, of $\mu$ associated to the
embedding $i$, the following identity in ${\bf R}/{\bf Z}$ holds,
\begin{align}\label{3.4}
H_X(\xi_{1,+},\xi_{1,-},V_1)=
H_X(\xi_{2,+},\xi_{2,-},V_2).\end{align}
\end{lemma}

{\it Proof.}  We first show that any direct image can be deformed
smoothly to another one which is of a simpler form.

 For any direct image
$(\xi_+,\xi_-,V=v+v^*)$   of $\mu$ associated to the embedding
$i:Y\hookrightarrow X$, let $N_\varepsilon$ be a tubular
neighborhood  of $Y$ which will be identified with a neighborhood
of the zero section in the total space of the normal bundle to $Y$
in $X$ for $\varepsilon$ small enough.

From the fundamental assumptions 1)-3) in Section \ref{2b}, over
$N_\varepsilon$ one has the identification
\begin{align}\label{3.54}
\xi_\pm|_{N_\varepsilon}\simeq\pi^*\left(\mu\otimes
S_\pm^*(N)\oplus F_\pm\right),
\end{align}
where $F_\pm=(\ker V|_Y)^{\bot}\cap \xi_\pm|_Y$.  Moreover, one
can write $v:\xi_+\to \xi_-$ near $Y$ in terms of $Z\in N$ such
that
\begin{align}\label{3.55}
v=\left(\begin{array}{cc}\pi^*{\rm Id}_\mu\otimes \tilde{c}(Z)+a & b\\
c & h\end{array}\right),\end{align} where $h:\pi^*F_+\to \pi^*F_-$
is an isomorphism, and the maps $b$ and $c$ have the infinitesimal
order $O(|Z|)$, while $a$ has the order $o(|Z|)$.

Choose a smooth cut-off function $f(Z)$ with support in
$N_{\varepsilon/2}$ and $f\equiv 1$ in $N_{\varepsilon/4}$. Set
for $0\leq s\leq 1,$
\begin{align}\label{3.56}
v_s=\left(\begin{array}{cc}\pi^*{\rm Id}_\mu\otimes \tilde{c}(Z)+a-sf(Z)a & b-sf(Z)b\\
c-sf(Z)c & h+sf(Z)(\pi^*(h|_Y)-h)\end{array}\right).
\end{align}

Clearly, for each $s\in [0,1]$, the map $v_s$ is globally
well-defined over $X$ and maps $\xi_+$ to $\xi_-$. Moreover, one
can choose $\varepsilon$ small enough so that $v_s$ is invertible
over $X\setminus Y$ for each $s$.

Note that $v_0=v$ and over $N_{\varepsilon/4}$
\begin{align}\label{3.57}
v_1=\left(\begin{array}{cc}\pi^*{\rm Id}_\mu\otimes \tilde{c}(Z) & 0\\
0 & \pi^*(h|_Y)\end{array}\right).
\end{align}

On the other hand, let $F_\pm$ carry Hermitian metrics $g^{F_\pm}$
and   Hermitian connections $\nabla^{F_\pm}$ respectively.

By using the identification (\ref{3.54}) one can choose another
Hermitian metric $g^\xi_1$ and Hermitian connection $\nabla^\xi_1$
on $\xi$ such that when restricted to $N_{\varepsilon/8}$, one has
that
\begin{align}\label{3.58}
g^{\xi_\pm}_1=\pi^*g^{\mu\otimes S^*_+(N)}\oplus
\pi^*g^{F_\pm},\quad \nabla^{\xi_\pm}_1=\pi^*\nabla^{\mu\otimes
S_+^*(N)}\oplus\pi^*\nabla^{F_\pm}.
\end{align}

For each $s\in [0,1]$, set
\begin{align}\label{3.60}
g^\xi_s=sg^\xi_1+(1-s)g^\xi,\quad \nabla^\xi_s={1\over
2}\left(\widetilde\nabla^\xi_s+\left(\widetilde\nabla^\xi_s\right)^*\right),\end{align}
where
\begin{align}\label{3.61}
\widetilde\nabla^\xi_s=s\nabla^\xi_1+(1-s)\nabla_0^\xi\end{align}
and $(\widetilde\nabla^\xi_s)^*$ is the adjoint of
$\widetilde\nabla^\xi_s$ with respect to the metric $g^\xi_s.$

Clearly for each $s\in [0,1]$,
\begin{align}\label{3.62}
\left(\xi,g^\xi_s,\nabla^\xi_s,V_s=v_s+v_s^*\right)\end{align} is
a direct image of $\mu$, where $v_s^*$ is the adjoint of $v_s$
with respect to $g^\xi_s$.

We now come back to the proof of the lemma.

From Lemma \ref{l3.2} and the above deformation, it is clear that
in order to prove the lemma at hand,  one needs only to prove it
for direct images $(\xi_{k,+},\xi_{k,-},V_k)$ with the properties
that
\begin{align}\label{3.5}
\xi_{k,\pm}|_{N_\varepsilon}\simeq \pi^*\left(
S^*_\pm(N)\otimes\mu\right)\oplus
\pi^*\left(F_{k,\pm}\right),\end{align}
\begin{align}\label{3.6}
g^{\xi_{k,\pm}|_{N_\varepsilon}}
=\pi^*\left.\left(g^{S_\pm^*(N)\otimes\mu}\oplus
g^{F_{k,\pm}}\right)\right|_{N_\varepsilon},
\nabla^{\xi_k,\pm|_{N_\varepsilon}}=\pi^*\left.\left(\nabla^{S^*_\pm(N)\otimes
\mu}\oplus
\nabla^{F_{k,\pm}}\right)\right|_{N_\varepsilon}\end{align}and
\begin{align}\label{3.7}
V_k|_{N_\varepsilon} =\tau^{N*}\widetilde{c}(Z)\oplus\pi^*{
h}_{F_k}\end{align} over $N_\varepsilon$ for $\varepsilon$ small
enough, where $h_{F_k}$ are self-adjoint odd automorphisms of
$F_k$, $k=1,\, 2$.

Now consider the super vector bundle
$$\left(\widetilde{\xi}_+,\widetilde{\xi}_-,\widetilde{V}\right)=\left(\xi_{1,+}\oplus
\xi_{2,-},\xi_{1,-}\oplus\xi_{2,+}, V_1\oplus V_2\right).$$ It is
easy to see that both $\widetilde{\xi}_\pm$ contain $\pi^*(
S^*(N)\otimes\mu)$.

 Clearly
over $N_\varepsilon$, one has
\begin{align}\label{3.8}
\widetilde{V}=\left(\begin{array}{cc}\pi^*{\rm Id}_\mu\otimes\tau^{N*}\tilde{c}(Z) & 0\\
0 & \pi^*{h}_{F_1 }\oplus\pi^*{h}_{F_2
}\end{array}\right).\end{align}

Define
\begin{align}\label{3.9}
i_g:\sqrt{-1}g(Z){\rm Id}_{\pi^*(\mu\otimes S^*(N))}\oplus
0_F:\widetilde{\xi}_+\to \widetilde{\xi}_-,\end{align} where $g$
is a cut-off function with support in $N_\varepsilon$ and $g\equiv
1$ in $N_{\varepsilon/2}$, and $0_F$ is the zero map from
$\pi^*(F_{1,+}\oplus F_{2,-})$ to $\pi^*(F_{1,-}\oplus F_{2,+})$.
Then  \begin{align}I_g=i_g+i_g^*\end{align} is a globally defined
odd endomorphism of $\widetilde{\xi}_+\oplus\widetilde{\xi}_-$.

Set \begin{align}\widetilde{V}_g=\widetilde{V}+I_g.\end{align} One
verifies that
\begin{align}\label{3.10}
\left. \left(\widetilde{V}_g\right)^2\right|_{N_\varepsilon}=\left(\begin{array}{cc}g(Z)^2+|Z|^2 & 0\\
0 & \pi^*{h}_{F_1 }\oplus\pi^*{h}_{F_2
}\end{array}\right).\end{align}

Note that $\widetilde{V}_g$ is invertible over $X$. By a result of
Bismut-Cheeger   (cf. \cite[Theorem 2.28]{BC}), one has
\begin{align}\label{3.11}
\overline{\eta}\left(D^{\widetilde{\xi}_+}\right)-\overline{\eta}\left(D^{\widetilde{\xi}_-}\right)\equiv
\int_X\widehat{A}\left(TX,\nabla^{TX}\right)\gamma^{\widetilde{\xi},\tilde{V}_g}\
\ {\rm mod}\ {\bf Z}.\end{align}

Let $F(s,T)$ be the curvature of the smooth family of the
superconnections
\begin{align}\label{3.12}
A(s,T)=\nabla^{\widetilde{\xi}}+\sqrt{T}\left(\widetilde{V}+sI_g\right)
\end{align}
on $\tilde{\xi}$.

For any closed form $\omega$ on $X$, by the standard double
transgression formula for Chern character forms (cf.
\cite[Proposition 3.1]{PV}), one has

\begin{multline}\label{3.13}
\int_X \omega\int_0^R\frac{1}{\sqrt{2\pi i}} \varphi{\rm
tr}_s\left[\widetilde{V}_ge^{-\left(\nabla^{\widetilde{\xi}}+\sqrt{T}\widetilde{V}_g\right)^2}\right]
\frac{dT}{2\sqrt{T}}\\
-\int_X\omega\int_0^R \frac{1}{\sqrt{2\pi i}} \varphi{\rm
tr}_s\left[\widetilde{V}e^{-\left(\nabla^{\widetilde{\xi}}+\sqrt{T}\widetilde{V}\right)^2}\right]
\frac{dT}{2\sqrt{T}}\\
=\int_X \omega \int_0^1 ds\frac{1}{\sqrt{2\pi i}}\varphi\, {\rm
tr}_s\left[\sqrt{R}I_ge^{-F(s,R)}\right].
\end{multline}

Note that when $R\to \infty$, the left hand side of the above equality converges to
\begin{align}\label{3.14}
\int_X\omega\gamma^{\widetilde{\xi},\widetilde{V}_g}-\int_X\omega\gamma^{\widetilde{\xi},\widetilde{V}}.
\end{align}

Now we want to show that the right hand side of (\ref{3.13}) tends
to zero as $R\to \infty$.

 Firstly, one sees that in $X\setminus N_{\varepsilon/2}$,
$$\int_0^1 ds\frac{1}{\sqrt{2\pi i}}\varphi\, {\rm
tr}_s\left[\sqrt{R}I_ge^{-F(s,R)}\right]$$
decays exponentially as $R\to \infty$.

On the other hand, since $g\equiv 1$ in $N_{\varepsilon/2}$, one
has
\begin{multline}
{\rm tr}_s\left[\sqrt{R}I_ge^{-F(s,R)}\right]=\sqrt{R}{\rm tr}_s\left[I_ge^{-(A(0,R)+\sqrt{R}sI_g)^2}\right]\\
=\sqrt{R}{\rm tr}_s\left[I_ge^{-A(0,R)^2}\right]e^{-Rs^2}.
\end{multline}

Since $A(0,R)^2$ maps $\xi_1$ to $\xi_1$ and $\xi_2$ to $\xi_2$
respectively, while $I_g$ exchange $\xi_1$ and $\xi_2$, one gets in $N_{\varepsilon/2}$ that
\begin{align}\label{3.15}
{\rm tr}_s\left[\sqrt{R}I_ge^{-F(s,R)}\right]\equiv 0.
\end{align}
Hence for any closed form $\omega$, one has
\begin{align}\label{3.16}
\int_X\omega\gamma^{\widetilde{\xi},\widetilde{V}_g}=\int_X\omega\gamma^{\widetilde{\xi},\widetilde{V}}.
\end{align}

By (\ref{3.11}) and (\ref{3.16}), one gets
\begin{align}\label{3.17}
\overline{\eta}\left(D^{\widetilde{\xi}_+}\right)-\overline{\eta}\left(D^{\widetilde{\xi}_-}\right)\equiv
\int_X\widehat{A}\left(TX,\nabla^{TX}\right)\gamma^{\widetilde{\xi},\widetilde{V}}\
\ {\rm mod}\ {\bf Z},
\end{align}which implies (\ref{3.4}).\ \ \ Q.E.D.

 $\ $

From the two lemmas above, one deduces easily that the
$H$-quantity of a direct image on $X$ of $\mu$ depends only on
$\mu$ as well as the isotropy class of the embedding. From now on,
we may denote this quantity by $H(Y,X,\mu)$, that is
\begin{align}\label{aa}
H_X(\xi_+,\xi_-,V)=H(Y,X,\mu).
\end{align}
Moreover, one sees that this quantity now does not depend on the
metrics and connections define it and is indeed a smooth
invariant!

\subsection{Direct images and Chern-Simons currents associated
with successive embeddings}\label{3a}

Let $(\xi_+,\xi_-,V)$ be a geometric direct image of $\mu$ for
$i:Y\hookrightarrow X$ constructed in \ref{2b}.

Let $j:X\hookrightarrow M$  be a totally geodesic embedding of $X$
into an odd-dimensional closed   oriented spin Riemannian manifold
$M$. Since the embedding $i:Y\hookrightarrow X$ is totally
geodesic, the induced embedding $j\circ i:Y\hookrightarrow M$ is
also totally geodesic.

Let $N_X$ be the normal bundle to $X$ in $M$. One can take
$\varepsilon_0>0$ appearing in Section \ref{2b} small enough so
that $N_{X,2\varepsilon_0}=\{u\in N_X: |u|<2 \varepsilon_0\}$  is
diffeomorphic to an open neighborhood of $X$ in $M$. Without
confusion we now view directly $N_{X,2\varepsilon_0}$ as an open
neighborhood of $X$ in $M$.

Let $N_Y$ be the normal bundle to $Y$ in $M$. Clearly,

\begin{align}\label{3.25}
N_Y=N\oplus i^*N_X.\end{align}
Then one can choose $\varepsilon_0>0$
small enough so that $N_{2\varepsilon_0}\oplus
i^*N_{X,2\varepsilon_0}$ is diffeomorphic to an open neighborhood of
$Y$ in $M$. Without confusion we now view directly
$N_{2\varepsilon_0}\oplus i^*N_{X,2\varepsilon_0}$ as an open
neighborhood of $Y$ in $M$.

Let $(\zeta_{+,+},\zeta_{+,-},W_+)$ (resp.
$(\zeta_{-,+},\zeta_{-,-},W_-)$) be the geometric direct image of
$\xi_+$ (resp. $\xi_-$) in the sense of Section \ref{2b}.

Let $\zeta=\zeta_+\oplus \zeta_-$ be the ${\bf Z}_2$-graded
Hermitian vector bundle over $M$ such that
\begin{align}\label{3.26}
\zeta_+= \zeta_{+,+}\oplus \zeta_{-,-},\ \ \ \zeta_-=
\zeta_{+,-}\oplus \zeta_{-,+}.\end{align} Then $\zeta_+-\zeta_-\in
\widetilde{K}(M)$ is a representative of $(j\circ i)_!\mu$ in the
sense of \cite{AH}.

Let
\begin{align}\label{3.27}
W=W_+\oplus W_-\end{align} be the induced odd endomorphism on
$\zeta$. Then ${\rm Supp}(W)=X$.

Let $f\in C^\infty(M)$ be such that ${\rm Supp}(f)\subset
N_{X,2\varepsilon_0}$ and $f\equiv 1$ on $N_{X,\varepsilon_0}$.

Let $\pi_X:N_X\rightarrow X$ denote the canonical projection.

Recall that by the construction of geometric direct images in
Section \ref{2b}, one has
\begin{align}\label{3.28}
\pi_X^*\left(S^*(N_X)\widehat{\otimes} \xi
\right)|_{N_{X,2\varepsilon_0}}\subset \zeta
|_{N_{X,2\varepsilon_0}}.
\end{align}

For any $Z\in N_{X,2\varepsilon_0}$, set

\begin{align}\label{3.29}
V_f(Z)=f(Z)\pi_X^*\left( {\rm
Id}_{S^*(N_X)} \widehat{\otimes}
 V\right):\pi_X^*\left(S^*(N_X)\widehat{\otimes}
\xi\right)\rightarrow\pi_X^*\left(S^*(N_X)\widehat{\otimes}
 \xi\right).\end{align}

 By (\ref{3.28}) it can be viewed  as an
endomorphism of $\zeta|_{N_{X,2\varepsilon_0}}$, which vanishes on
the orthogonal  complement of $\pi_X^* (S^*(N_X)\widehat{\otimes}
 \xi  )$ in $\zeta|_{N_{X,2\varepsilon_0}}$.

 Also, since ${\rm
Supp}(f)\subset N_{X,2\varepsilon_0}$, one can extend $V_f $ as
zero endomorphism to $M\setminus N_{X,2\varepsilon_0}$.

Thus, we may view $V_f $ as an odd self-adjoint endomorphism on
$\zeta$.

Let $W_f$ be the odd self adjoint endomorphism of $\zeta$ defined
by
\begin{align}\label{3.30}
W_f=W+V_f.
\end{align}
By (\ref{2.4}), (\ref{3.27}), (\ref{3.29}) and (\ref{3.30}), one
verifies easily that
\begin{align}\label{3.31}
{\rm Supp}\left(W_f\right)= Y.
\end{align}
Moreover, by (\ref{3.25}), one sees that the obvious analogues of
the conditions \cite[(1.10) and (1.12)]{BZ} verify for
$(\zeta_+,\zeta_-,W_f)$ near $Y\subset M$.

Thus, one gets a well-defined Chern-Simons current $\gamma^{\zeta,
W_f}$ over $M$ as in (\ref{2.16}).

\begin{prop}\label{t3.3} If $g$ is another cut-off function
verifying the same condition as $f$, then up to smooth exact forms
on $M$, $\gamma^{\zeta, W_f}=\gamma^{\zeta, W_g}$.\end{prop}

{\it Proof}. Take $f_s=f+s(g-f)$, $0\leq s\leq 1$. Then $f_s$ is a
family of cut-off functions verifying the same condition as $f$.
Let
\begin{align}\label{3.32}
A(s,T)=\nabla^\zeta+\sqrt{T}W_{f_s}
\end{align} be the corresponding smooth family of superconnections
on $\zeta$. Let $F(s,T)$ be the curvature of $A(s,T)$.

By the double transgression formula for Chern character forms (cf.
\cite[Proposition 3.1]{PV}), one has
\begin{align}\label{3.33}
\frac{d}{ds} \mathrm{tr}_s\left[\frac{dA(s,T)}{dT}e^{-{F(s,T)}}\right]-\frac{d}{dT}
\mathrm{tr}_s\left[ \frac{dA(s,T)}{ds}e^{- {F(s,T)}}\right]
=dv(s,T),\end{align}
where
\begin{multline}
v(s,T)=-\int_{0}^{1} \mathrm{tr}_s\left[\frac{dA(s,T)}{ds}
e^{-u {F}(s,T)}\frac{d A(s,T)}{dT} e^{-(1-u) {F(s,T)}}\right]du \\
=- \int_{0}^{1} \mathrm{tr}_s\left[\left(\sqrt{T}(g -f ){\rm
Id}_{\pi_X^*S^*(N_X)}\widehat{\otimes}\pi_X^*(V)
\right)e^{-u {F(s,T)}}\frac{dA(s,T)}{dT} e^{-(1-u) {F(s,T)}}du\right] \\
=-(g-f)\int_{0}^{1} \mathrm{tr}_s\left[\left(\sqrt{T} {\rm
Id}_{\pi_X^*S^*(N_X)}\widehat{\otimes}\pi_X^*(V) \right)e^{-u
{F(s,T)}}\frac{dA(s,T)}{dT} e^{-(1-u) {F(s,T)}}du\right].
\end{multline}

Since $$\int_0^1ds\int_0^{+\infty}dT \int_{0}^{1}\mathrm{tr}_s\left[ \sqrt{T}\pi_X^*
\left({\rm Id}_{S^*(N_X)}\widehat{\otimes} V\right)e^{-u{F(s,T)}}\frac{dA(s,T)}{dT} e^{-(1-u) {F(s,T)}}du\right]$$
is a smooth form on $N_{X,\varepsilon_0} \setminus X$ (cf. \cite[Proposition 3.8]{PV}) and $g-f $ vanishes near $X$, one
sees that
$$\int_0^1ds\int_0^{+\infty}dTv(s,T)$$
is a smooth form. Thus, one has
\begin{multline}
\gamma^{\zeta,W_{g}}-\gamma^{\zeta,W_{f}}=\lim_{R\to+\infty}\int_0^1ds
\frac{d}{ds}\left(\int_0^R\frac{1}{\sqrt{2\pi i}}\varphi\,
\mathrm{tr}_s\left[\frac{dA(s,T)}{dT} e^{- {F(s,T)}}\right]dT\right)\\
=\lim_{R\to+\infty}\int_0^1ds\int_0^R\frac{1}{\sqrt{2\pi
i}}\frac{d}{dT}\varphi\,\mathrm{tr}_s\left[\frac{dA(s,T)}{ds} e^{-
{F}(s,T)}\right]dT\\
+\frac{1}{\sqrt{2\pi i}}\varphi\int_0^1ds\int_0^{+\infty}dv(s,T)\\
=\int_0^1ds \lim_{T\to+\infty}\frac{1}{\sqrt{2\pi
i}}\varphi\,\mathrm{tr}_s\left[\left(\sqrt{T}(g-f )\pi_X^*
\left({\rm Id}_{S^*(N_X)}\widehat{\otimes}
V\right)\right)e^{-{F(s,T)}}\right]\\
+ \frac{1}{\sqrt{2\pi i}}\varphi\, d\int_0^1ds\int_0^{+\infty}v(s,T)\\
=\frac{1}{\sqrt{2\pi i}}\varphi\,
d\int_0^1ds\int_0^{+\infty}v(s,T),
\end{multline}
from which Proposition \ref{t3.3} follows. \ \ Q.E.D.

\subsection{Real embeddings and Chern-Simons current: a
Riemann-Roch formula}\label{3b}

In this section, we prove the following result, which might be
thought of as a Riemann-Roch type formula for the Chern-Simons
currents $\gamma^{\zeta,W_f}$ and $\gamma^{\zeta,
W}=\gamma^{\zeta_+,W_+}-\gamma^{\zeta_-,W_-}$.

\begin{thm}\label{t3.4}  For any closed form $\omega$ on $M$, one
has
\begin{align}\label{3.34}
\int_M\omega\gamma^{\zeta,W_f}-\int_M\omega\gamma^{\zeta,W}=
\int_X\left(j^*\omega\right)\widehat{A}^{-1}(N_X)\gamma^{\xi,V}.\end{align}
In other words,
$\gamma^{\zeta,W_f}-\gamma^{\zeta,W}-\widehat{A}^{-1}(N_X)\gamma^{\xi,V}\delta_X$
is a current on $M$ eliminating closed forms.
\end{thm}

{\it Proof}. By using formula (\ref{3.33}) for $f_s=sf$, one sees
that for any closed form $\omega$ on $M$, one has
\begin{multline}\label{3.35}
\int_M\omega\int_0^R\frac{1}{\sqrt{2\pi
i}}\varphi\,\mathrm{tr}_s\left[ {W_f} e^{-
{\left(\nabla^\zeta+\sqrt{T}W_f\right)^2}}\right]{dT\over
2\sqrt{T}}\\ -\int_M\omega\int_0^R\frac{1}{\sqrt{2\pi
i}}\varphi\,\mathrm{tr}_s\left[ {W }
e^{- {\left(\nabla^\zeta+\sqrt{T}W \right)^2}}\right]{dT\over 2\sqrt{T}}\\
= \int_M\omega\int_0^1ds  \frac{1}{\sqrt{2\pi
i}}\varphi\,\mathrm{tr}_s\left[\left(\sqrt{R} f  \pi_X^*\left({\rm
Id}_{S^*(N_X)}\widehat{\otimes} V\right)\right)e^{-
{F(s,R)}}\right]  .\end{multline}

It is clear that in order to prove (\ref{3.34}), one needs to
evaluate the limit as $R\rightarrow +\infty$ of the right hand side
of (\ref{3.35}).

In the current case
\begin{align}\label{3.36}
A(s,T)=\nabla^\zeta+\sqrt{T}W_{sf
}=\nabla^\zeta+\sqrt{T}\left(W+sV_f\right).
\end{align}

From (\ref{3.29}),  (\ref{3.36}), the construction of
$(\zeta_+,\zeta_-,W)$ and the fact that $f\equiv 1$ on
$N_{X,\varepsilon_0}$, one deduces that for any $T\geq 0$, $0\leq
s\leq 1$, one has on $N_{X,\varepsilon_0}$ that
\begin{multline}\label{3.37}
\mathrm{tr}_s\left[\left(\sqrt{T} f  \pi_X^*\left({\rm
Id}_{S^*(N_X)}\widehat{\otimes} V\right)\right)e^{-
{F(s,T)}}\right]\\
=\pi_X^*\mathrm{tr}_s\left[\sqrt{T}Ve^{-\left(\nabla^\xi+s\sqrt{T}V\right)^2}\right]\cdot
\mathrm{tr}_s\left[
e^{-\left(\nabla^{\pi^*_XS(N_X)}+\sqrt{T}\tau^{
{N_X}*}\widetilde{c}(Z)\right)^2}\right].
\end{multline}

Let $\psi\geq 0$ be a smooth function on $M$ such that ${\rm
Supp}(\psi)\subset N_{X,\varepsilon_0}$, $\psi\equiv 1$ on
$N_{X,{1\over 2}\varepsilon_0}$.

By (\ref{3.37}), one has
\begin{multline}\label{3.38}
 \int_M\omega\int_0^1ds  \frac{1}{\sqrt{2\pi
i}}\varphi\,\mathrm{tr}_s\left[\left(\sqrt{T} f  \pi_X^*\left({\rm
Id}_{S^*(N_X)}\widehat{\otimes} V\right)\right)e^{-
{F(s,T)}}\right]  \\
= \int_{N_{X,\varepsilon_0}}\psi\omega\int_0^1ds
\frac{1}{\sqrt{2\pi i}}\varphi\,\mathrm{tr}_s\left[\left(\sqrt{T}
f \pi_X^*\left({\rm Id}_{S^*(N_X)}\widehat{\otimes}
V\right)\right)e^{- {F(s,T)}}\right] \\ +
\int_M(1-\psi)\omega\int_0^1ds \frac{1}{\sqrt{2\pi
i}}\varphi\,\mathrm{tr}_s\left[\left(\sqrt{T} f \pi_X^*\left({\rm
Id}_{S^*(N_X)}\widehat{\otimes} V\right)\right)e^{-
{F(s,T)}}\right]   ,
\end{multline}
with
\begin{multline}\label{3.39}
 \int_{N_{X,\varepsilon_0}}\psi\omega\int_0^1ds  \frac{1}{\sqrt{2\pi
i}}\varphi\,\mathrm{tr}_s\left[\left(\sqrt{T} f  \pi_X^*\left({\rm
Id}_{S^*(N_X)}\widehat{\otimes} V\right)\right)e^{-
{F(s,T)}}\right]\\
=\int_X \int_0^1\frac{ds}{\sqrt{2\pi
i}}\varphi\,\mathrm{tr}_s\left[\sqrt{T}Ve^{-\left(\nabla^\xi+s\sqrt{T}V\right)^2}\right]
\int_{N_X/X}\psi\omega\varphi\,\mathrm{tr}_s\left[
e^{-\left(\nabla^{\pi^*_XS(N_X)}+\sqrt{T}\tau^{
{N_X}*}\widetilde{c}(Z)\right)^2}\right]\\
=\int_X \int_0^T {1\over\sqrt{2\pi i}}\varphi\,\mathrm{tr}_s\left[
Ve^{-\left(\nabla^\xi+\sqrt{s}V\right)^2}\right]{ds\over
2\sqrt{s}}\int_{N_X/X}\psi\omega\varphi\,\mathrm{tr}_s\left[
e^{-\left(\nabla^{\pi^*_XS(N_X)}+\sqrt{T}\tau^{
{N_X}*}\widetilde{c}(Z)\right)^2}\right],
\end{multline}
where $\int_{N_X/X}$ is the integration of differential forms along
the fibre of $N_X$ over $X$.

By proceeding as in \cite[Theorem 1.2]{BZ}, one sees that there
exists $C>0$ such that as $T>0$ is large enough,
\begin{multline}\label{3.40}
\left\|\int_{N_X/X}\psi\omega\varphi\mathrm{tr}_s\left[
e^{-\left(\nabla^{\pi^*_XS(N_X)}+\sqrt{T}\tau^{
{N_X}*}\widetilde{c}(Z)\right)^2}\right]-\left(j^*\omega\right)
\widehat{A}^{-1}(N_X)\right\|_{C^1(X)}\\ \leq {C\over
\sqrt{T}}\|\omega\|_{C^2(M)}.
\end{multline}

By (\ref{3.40}) and \cite[Theorem 1.2]{BZ} again, one sees that as
$T\rightarrow +\infty$,
\begin{multline}\label{3.41}
\left|\int_X
\left(\int_{N_X/X}\psi\omega\varphi\,\mathrm{tr}_s\left[
e^{-\left(\nabla^{\pi^*_XS(N_X)}+\sqrt{T}\tau^{
{N_X}*}\widetilde{c}(Z)\right)^2}\right]-\left(j^*\omega\right)
\widehat{A}^{-1}(N_X)\right)\right.\\ \cdot\left.\int_0^{T}
{1\over\sqrt{2\pi i}}\varphi\,\mathrm{tr}_s\left[
Ve^{-\left(\nabla^\xi+\sqrt{s}V\right)^2}\right]{ds\over
2\sqrt{s}}\right|\rightarrow 0.
\end{multline}

From (\ref{2.12}), (\ref{3.39}) and (\ref{3.41}), one finds that as
$T\rightarrow +\infty$,
\begin{multline}\label{3.42}
 \int_{N_{X,\varepsilon_0}}\psi\omega\int_0^1ds  \frac{1}{\sqrt{2\pi
i}}\varphi\,\mathrm{tr}_s\left[\left(\sqrt{T} f  \pi_X^*\left({\rm
Id}_{S^*(N_X)}\widehat{\otimes} V\right)\right)e^{-
{F(s,T)}}\right]\\ \rightarrow \int_X
(j^*\omega)\widehat{A}^{-1}(N_X)\gamma^{\xi,V}.
\end{multline}

On the other hand, if we write
\begin{multline}\label{3.43}
 \int_M(1-\psi)\omega\int_0^1ds \frac{1}{\sqrt{2\pi
i}}\varphi\,\mathrm{tr}_s\left[\left(\sqrt{T} f \pi_X^*\left({\rm
Id}_{S^*(N_X)}\widehat{\otimes} V\right)\right)e^{-
{F(s,T)}}\right] \\ =  \int_X \int_0^T {1\over\sqrt{2\pi
i}}\varphi\,\mathrm{tr}_s\left[
Ve^{-\left(\nabla^\xi+\sqrt{s}V\right)^2}\right]{ds\over
2\sqrt{s}}\\
\cdot\int_{N_X/X}(1-\psi)f\omega\varphi\,\mathrm{tr}_s\left[
e^{-\left(\nabla^{\pi^*_XS(N_X)}+\sqrt{T}\tau^{
{N_X}*}\widetilde{c}(Z)\right)^2}\right],
\end{multline}
then by noting that as $T\rightarrow +\infty$, $$\int_0^T
{1\over\sqrt{2\pi i}}\varphi\,\mathrm{tr}_s\left[
Ve^{-\left(\nabla^\xi+\sqrt{s}V\right)^2}\right]{ds\over
2\sqrt{s}}
$$
grows polynomially in $T$ (compare with \cite[Lemma 6.1]{PV}),
while since $1-\psi\equiv 0$ on $N_{X,{1\over 2}\varepsilon}$,
$$\int_{N_X/X}(1-\psi)f\omega\varphi\,\mathrm{tr}_s\left[
e^{-\left(\nabla^{\pi^*_XS(N_X)}+\sqrt{T}\tau^{
{N_X}*}\widetilde{c}(Z)\right)^2}\right]$$ decays exponentially in
$T$, one sees that as $T\rightarrow +\infty$,
\begin{align}\label{3.44}
 \int_M(1-\psi)\omega\int_0^1ds \frac{1}{\sqrt{2\pi
i}}\varphi\,\mathrm{tr}_s\left[\left(\sqrt{T} f \pi_X^*\left({\rm
Id}_{S^*(N_X)}\widehat{\otimes} V\right)\right)e^{-
{F(s,T)}}\right]\rightarrow 0.
\end{align}

From (\ref{3.35}), (\ref{3.38}), (\ref{3.42}) and (\ref{3.44}), one
gets (\ref{3.34}), which completes the proof of Theorem \ref{t3.4}.\
\ Q.E.D.

$\ $

\begin{Rem}\label{t3.5} From (\ref{2.13}) and (\ref{3.34}), one formally gets
\begin{multline}\label{3.45}
d\gamma^{\zeta,W_f}={\rm ch}\left(\zeta,\nabla^\zeta\right)
+\widehat{A}^{-1}\left(N_X,\nabla^{N_X}\right)\left(d\gamma^{\xi,V}-{\rm
ch}\left(\xi,\nabla^\xi\right)\right)\delta_X\\ ={\rm ch}\left(\zeta,\nabla^\zeta\right)
-\widehat{A}^{-1}\left(N_Y,\nabla^{N_Y}\right){\rm ch}\left(\mu,\nabla^\mu\right) \delta_Y,
\end{multline} which fits with (\ref{2.13}) again.
\end{Rem}

$\ $

\begin{cor}\label{c3.1} Under the assumptions and notations above, the following
identity in ${\bf R}/{\bf Z}$ holds,
\begin{align}\label{7.1}
H(Y,M,\mu)=H(Y,X,\mu)+H(X,M,\xi_+) -H(X,M,\xi_-).
\end{align}
\end{cor}

{\it Proof}. By (\ref{3.1})  and (\ref{3.34}), one has
\begin{multline}\label{7.2}
H_M\left(\zeta_+,\zeta_-,W_f\right)\equiv
\overline{\eta}\left(D^{\zeta_+}\right)-\overline{\eta}\left(D^{\zeta_-}\right)
-\int_M\widehat{A}\left(TM,\nabla^{TM}\right)\gamma^{\zeta,W_f}-
\overline{\eta}\left(D^\mu\right)\ \ {\rm mod}\ {\bf Z}\\
\equiv
\overline{\eta}\left(D^{\zeta_{++}}\right)+\overline{\eta}\left(D^{\zeta_{+-}}\right)
-\overline{\eta}\left(D^{\zeta_{-+}}\right)
-\overline{\eta}\left(D^{\zeta_{+-}}\right)-\int_M\widehat{A}\left(TM,\nabla^{TM}\right)\gamma^{\zeta,W}\\
- \int_X\widehat{A}\left(TX,\nabla^{TX}\right)\gamma^{\xi,V}-
\overline{\eta}\left(D^\mu\right)\ \ {\rm mod}\ {\bf Z}\\
\equiv
\overline{\eta}\left(D^{\zeta_{++}}\right)-\overline{\eta}\left(D^{\zeta_{+-}}\right)-
\int_M\widehat{A}\left(TM,\nabla^{TM}\right)\gamma^{\zeta_+,W_+}-\overline{\eta}\left(D^{\xi_+}\right)\\
-\left(\overline{\eta}\left(D^{\zeta_{-+}}\right)-\overline{\eta}\left(D^{\zeta_{--}}\right)-
\int_M\widehat{A}\left(TM,\nabla^{TM}\right)\gamma^{\zeta_-,W_-}-\overline{\eta}\left(D^{\xi_-}\right)\right)\\
+\overline{\eta}\left(D^{\xi_{+}}\right)-\overline{\eta}\left(D^{\xi_{-}}\right)
- \int_X\widehat{A}(TX,\nabla^{TX})\gamma^{\xi,V}-
\overline{\eta}\left(D^\mu\right)\ \ {\rm mod}\ {\bf Z}\\
=H_M\left(\zeta_{++},\zeta_{+-},W_+\right)
-H_M\left(\zeta_{-+},\zeta_{--},W_-\right)+H_X\left(\xi_+,\xi_-,V\right).
\end{multline}

By (\ref{aa}) and (\ref{7.2}), one gets (\ref{7.1}), which
completes the proof of Corollary \ref{c3.1}. \ \ Q.E.D.

\subsection{Proof of the Bismut-Zhang localization formula}\label{3c}

 Recall that by Remark \ref{t2.3} one knows that there exists
a totally geodesic embedding $i_0:Y\hookrightarrow S^{2m-1}$  such
that
\begin{align}\label{6.1}
H \left(Y,S^{2m-1},\mu\right)=0.\end{align}

We first show that this holds for any embedding of $Y$ to an
arbitrary  sphere.

\begin{lemma}\label{t6.1} Let $\mu$ be a  complex vector bundle  over an odd dimensional
closed oriented spin  manifold $Y$.  Then for any embedding
$i:Y\hookrightarrow S^{2n-1}$, the following identity in ${\bf
R}/{\bf Z}$ holds,
\begin{align}\label{6.2}
H \left(Y,S^{2n-1},\mu\right)=0.\end{align}
\end{lemma}

{\it Proof}. For the two embeddings $i_0:Y\hookrightarrow
S^{2m-1}$ and $i:Y\hookrightarrow S^{2n-1}$, we first consider the
associated embeddings $i_0':Y\hookrightarrow S^{2m-1}(1)\subset
{\bf R}^{2m}$ and $i':Y\hookrightarrow S^{2n-1}(1)\subset {\bf
R}^{2n}$ to the standard unit  spheres.

By using a trick in  \cite[Page 498]{AS}, we construct a smooth
family of embeddings  $j_s:Y\hookrightarrow S^{2m+2n-1}(1)$,
$0\leq s\leq 1$, obtained by
\begin{align}\label{6.3}
y\in Y\mapsto {1\over \sqrt{s^2+(1-s)^2}}\left((1-s)i_0'y,
si'y\right)\in S^{2m+2n-1}(1)\subset {\bf R}^{2m}\oplus{\bf
R}^{2n}.
\end{align}

Now since $j_0$ and $j_1$ are isotropic to each other by
(\ref{6.3}), by Lemma \ref{l3.2} one has
\begin{align}\label{6.4}
H_{j_0} \left(Y,S^{2m+2n-1},\mu\right)=H_{j_1}
\left(Y,S^{2m+2n-1},\mu\right).
\end{align}

On the other hand, by the Bott periodicity, any complex vector
bundle over an odd dimensional sphere can be expressed as a
difference of two trivial vector bundles so that can be extended
to a bounding ball, one can then apply the arguments in \cite{B1}
and \cite{Z2} to see that
\begin{align}\label{6.5}
H  \left(S^{2k-1},S^{2k'-1},\nu\right)=0
\end{align}
for any standard embedding between spheres and $\nu$ any complex
vector bundle over $S^{2k-1}$.

Now by applying Corollary \ref{c3.1} and (\ref{6.5}) to the
successive embedding $j_0 :Y\hookrightarrow
S^{2m-1}\hookrightarrow S^{2m+2n-1}$ and  $j_1 :Y\hookrightarrow
S^{2n-1}\hookrightarrow S^{2m+2n-1}$ respectively, one gets
\begin{align}\label{6.6}
H_{j_0 }
\left(Y,S^{2m+2n-1},\mu\right)=H\left(Y,S^{2m-1},\mu\right)=0
\end{align}
and
\begin{align}\label{6.7}
H_{j_1 } \left
(Y,S^{2m+2n-1},\mu\right)=H\left(Y,S^{2n-1},\mu\right)
\end{align} respectively.

From (\ref{6.4}),  (\ref{6.6}) and  (\ref{6.7}), one gets
(\ref{6.2}). \ \ Q.E.D.

$\ $

We now come to the proof of the Bismut-Zhang localization Theorem
\ref{t2.1} which is equivalent to saying that for any embedding
$i:Y\hookrightarrow X$ between two odd dimensional closed oriented
spin
  manifolds and a complex vector bundle $\mu$ over $Y$,
one has
\begin{align}\label{6.8}
H   (Y,X,\mu)=0.
\end{align}

Indeed, let $j:X\hookrightarrow S^{2N-1}$ be a further embedding
of $X$ into a higher odd dimensional sphere, let
$\xi=\xi_+\oplus\xi_-$ be a ${\bf Z}_2$-graded vector bundle over
$X$ so that by giving suitable metrics, connections and odd
endomorphism $V$ of $\xi$, $(\xi_+,\xi_-,V)$ realizes  a direct
image of $\mu$ in the sense of Section \ref{2b}.

By Corollary \ref{c3.1} and Lemma \ref{t6.1}, one deduces that
\begin{align}\label{6.9}
H
(Y,X,\mu)=H\left(Y,S^{2N-1},\mu\right)-H\left(X,S^{2N-1},\xi_+\right)+H\left(X,S^{2N-1},\xi_-\right)=0,
\end{align}
which via (\ref{6.8}) completes the proof of Theorem \ref{t2.1}. \
\ Q.E.D.

\end{document}